\documentstyle[12pt]{article}
  
\newcommand{\numero}[1]{
\addtocounter{section}{1}
\begin{center}{\bf \thesection .\
#1\vspace{-.1in}}\nopagebreak\end{center}
\setcounter{subsection}{0}
\setcounter{lemma}{0}\indent}

\newcommand{\subnumero}[1]{
\pagebreak[1]\begin{center}{\em #1}\nopagebreak\end{center}
}

\newcommand{\eop}{\hfill $\Box$\vspace{.1in}}

\newtheorem{lemma}{Lemma}[section]
\newtheorem{theorem}[lemma]{Theorem}

\newtheorem{corollary}[lemma]{Corollary}

\newtheorem{proposition}[lemma]{Proposition}

\newcommand{\cc}{{\bf C}}

\newcommand{\zz}{{\bf Z}}

\newcommand{\Cc}{{\cal C}}
\newcommand{\Ee}{{\cal E}}

\newcommand{\Gg}{{\cal G}}

\newcommand{\Oo}{{\cal O}}
\newcommand{\Dd}{{\cal D}}

\newcommand{\Bb}{{\cal B}}

\newcommand{\Hh}{{\cal H}}

\newcommand{\Xx}{{\cal X}}
\newcommand{\Yy}{{\cal Y}}
\newcommand{\Zz}{{\cal Z}}

\begin{document}

\begin{center}
{\bf \Large
The topological realization of a simplicial presheaf}
\\
Carlos Simpson
\end{center}

\begin{center}
{\bf Introduction}
\end{center}

While preparing for EUROPROJ's Trento school on stacks (September
1996) it became apparent that  an obvious thing that one
would like to say about a stack---to take its topological
realization---was not altogether obvious to define or to handle.
The same question can be posed more generally for a presheaf of
spaces or (as it is common to say) a simplicial presheaf.

We give
a definition of the topological realization of a
presheaf of spaces on a site $\Xx$,
with respect to a covariant functor $F: \Xx \rightarrow Top$. It
is a topological space defined in a fairly obvious way. In our usual
case $\Xx$ will a site of schemes over $Spec (\cc )$ and $F$ is the
functor which to a scheme $X$ associates its underlying usual
topological space $X^{\rm top}$.  If $G$ is representable by an
object $X$ then the realization of $G$ is (homotopic to) $X^{\rm
top}$.

One recovers the topological realization of a simplicial presheaf
by first realizing over each object into a presheaf of spaces, and
then applying our definition.  One recovers the topological
realization of a stack (of groupoids) by first strictifying into a
presheaf of groupoids then taking the realization of the nerve of
the groupoid over each object to get a presheaf of spaces.

Once the theory of $n$-stacks is well off the ground, the same
remark will hold for the topological realization of an
$n$-stack (whereas for now one must replace the notion of
$n$-stack by $n$-truncated presheaf of spaces, and from this point
of view one can directly take the topological realization as
defined below).

The main theorem is the invariance of topological realization
under Illusie weak equivalence. Recall that if $G$ is a presheaf
of spaces on the site $\Xx$ then we obtain a presheaf of sets $\pi
_0^{\rm pre}(G)$ on $\Xx$ and for any $g\in G(X)$ a presheaf of
groups $\pi _i^{\rm pre}(G|_{\Xx /X}, g)$ on $\Xx /X$.
Then $\pi _0(G)$ (resp. $\pi _i(G|_{\Xx /X},g)$) is
defined to be the
sheafification  of  $\pi
_0^{\rm pre}(G)$ (resp. the sheafification of $\pi _i^{\rm
pre}(G|_{\Xx /X}, g)$ on $\Xx /X$).  A morphism $\psi :G\rightarrow
G'$ of presheaves of spaces is an {\em Illusie weak equivalence} if
it induces an isomorphism on the homotopy sheaves
$$
\pi _0(G)\stackrel{\cong}{\rightarrow}\pi _0(G')
$$
and for any $X\in \Xx$ and $g\in G(X)$ an isomorphism
$$
\pi _i^{\rm
pre}(G|_{\Xx /X}, g)\stackrel{\cong}{\rightarrow}\pi _i^{\rm
pre}(G'|_{\Xx /X}, \psi (g)).
$$
Heuristically, the correct homotopy type of presheaves of spaces on
$\Xx$ is the homotopy type up to Illusie weak equivalence (this is
how the Grothendieck topology is taken into account). More
concretely in the case of stacks for example, one often has a
presentation of a stack as the quotient of a scheme by a relation.
Taking the resulting simplicial scheme and then realizing over each
object to get a  presheaf of spaces gives something which is Illusie
equivalent to the presheaf  of spaces associated to the stack.
Thus we would like to define the realization in terms of a
presheaf of spaces but we would like to know that it is an
invariant of the homotopy type, i.e. we would like to know that
the morphism of realizations induced by an Illusie weak
equivalence is a weak homotopy equivalence of spaces---this is the
result of Theorem \ref{main}.

We give several examples of calculations of topological
realizations. These generally make use of Theorem \ref{main}.
The calculation that I originally wanted to look at for the Trento
school is that of the realization of a one-dimensional
Delinge-Mumford stack.

At the end we sketch how to extend this theory to the case of functors with
values in an $n$-topos (although we don't yet know very clearly what an
$n$-topos is).  This section resulted from email conversations with K. Behrend
and C. Telemann, who were both talking about topological realizations
but also pullbacks of stacks from one site to another.

{\em Caution:}  At the time of writing this first version, I don't
know if the main theorem has already appeared somewhere else.   In
the case of presheaves of sets over a topological space, I think
that it is already known as the statement that the ``espace etal\'e''
of a presheaf of sets is the same as the
``espace etal\'e'' of the associated sheaf.
I imagine that our main theorem has already been {\em used} at least
in the case of stacks, because it is such an obvious statement that
it is easy to use it without thinking to prove it (in particular the
reader will probably find that the calculations given in most of our
examples are nothing new). So if anybody knows of a reference where
the theorem (or some version of it) has already been proved, please
let me know (by email: carlos@picard.ups-tlse.fr)!

Let $Top$ denote the category of topological spaces. Let $\Delta ^n$
denote the standard $n$-simplex. In general $\ast$ denotes the
one-point space, or a constant functor to $Top$ whose values are the
one-point space.

\numero{The definition of topological realization}

Suppose $\Xx$ is a category, and suppose $F: \Xx \rightarrow Top$
is a covariant functor and $G: \Xx \rightarrow Top$ a
contravariant functor. Then we define the {\em realization of $G$ with
respect to $F$} denoted by $\Re _{\Xx}(F, G)$ as the quotient
$$
\Re _{\Xx}(F, G)= \coprod _{\phi _{\cdot} : X_0 \rightarrow X_n}
F(X_0)\times G(X_n)\times \Delta ^n / \sim
$$
by the equivalence relation $\sim $ to be explained below. The
notation $\phi _{\cdot} : X_0\rightarrow X_n$ means a composable
sequence of $n$-morphisms
$$
X_0\stackrel{\phi _1}{\rightarrow} X_1\ldots X_{n-1}\stackrel{\phi
_n}{\rightarrow}X_n
$$
in $X$.

The equivalence relation $\sim$ (the obvious one in our
situation) is defined by saying that  if $t\in \Delta ^m$ is a point and
if $\eta : \Delta ^m\rightarrow \Delta ^n$ is a face, then
$(a,b,\eta (t)) \in F(X_0)\times
G(X_n)\times \Delta ^n$ (indexed by $\phi _{\cdot}$) is equivalent to
$(a', b', t)\in F(X_i)\times G(X_{j})\times \Delta ^m$
(indexed by $\phi '_{\cdot}$), where $i$ and
$j$ are the numbers of the first and last vertices of the face $\eta$,
$\phi '_{\cdot}: X_i\rightarrow X_j$ is the composable sequence
obtained from $\phi _{\cdot}$ by applying the face $\eta$, and
$a'$ is the image  of $a$ under $F(X_0\rightarrow
F(X_i)$ and $b'$ is the image of $b$ under $G(X_n)\rightarrow
G(X_j)$.

The space $\Re _{\Xx}(F,G)$ is topologized as the quotient of the
disjoint union of the spaces $F(X_0)\times G(X_n)\times \Delta ^n$.
Note that if we order the addition of ``cells'' according to $n$, then
a ``cell'' of the form $F(X_0)\times G(X_n)\times \Delta ^n$ is  added
onto the previous part via an attaching map from $F(X_0)\times
G(X_n)\times (\partial \Delta ^n)$ to the previous part; the attaching
map is defined by the equivalence relation $\sim$.

A morphism $F\rightarrow F'$ (resp. $G\rightarrow G'$) is said to be
an {\em object-by-object equivalence} if it induces a weak homotopy
equivalence $F(X)\cong  F'(X)$ (resp. $G(X)\cong G'(X)$) for each
$X\in \Xx$.  If $F\rightarrow F'$ and $G\rightarrow G'$ are
object-by-object equivalences then they induce a weak homotopy
equivalence $\Re _{\Xx}(F,G)\rightarrow \Re _{\Xx}(F',G')$.

If $\pi : \Xx \rightarrow \Yy$ is a functor and if $F, G: \Yy
\rightarrow Top$ are as above, then we obtain a natural
transformation
$$
\Re _{\Xx }(\pi ^{\ast}(F), \pi ^{\ast}(G)) \rightarrow \Re
_{\Yy}(F,G).
$$
This is compatible with compositions $\Xx \rightarrow \Yy \rightarrow
\Zz$ in the obvious way.

\begin{lemma}
\label{finalObject}
Suppose $\Xx$ is a category with final object $U$. Suppose $F:\Xx
\rightarrow Top$ is a functor. Then the natural morphism
$$
F(U)\rightarrow \Re _{\Xx }(F, \ast )
$$
is a weak homotopy equivalence.
\end{lemma}
{\em Proof:}
In $\Re _{\Xx }(F, \ast )$ we homotope a cell of the form
$F(X_0)\times \Delta ^n$ corresponding to a composable $n$-tuple
$\phi _{\cdot}: X_0 \rightarrow X_n$ to $F(U)$ through the cell
$F(X_0)\times \Delta ^{n+1}$ corresponding to $(\phi _{\cdot},
\beta )$ where $\beta : X_n\rightarrow U$ is the canonical morphism.
Note that the last vertex of the cell $F(X_0)\times \Delta ^{n+1}$
gets attached to $F(U)$.  This homotopy is a retraction from
$\Re _{\Xx }(F, \ast )$ to $F(U)$.
\eop

\subnumero{Relative realizations}

If $\pi : \Xx \rightarrow \Yy$ is a functor (when nothing is specified
that means covariant!) and if $F, G: \Xx \rightarrow Top$ are as above,
then we can define the {\em standard relative realization} which is
a covariant functor
$$
\Re _{\Xx /\Yy}(F, G): \Yy \rightarrow Top.
$$
It is defined by
$$
\Re _{\Xx /\Yy}(F,G)(Y):= \Re _{\Xx /Y}(F|_{\Xx /Y},
G|_{\Xx /Y}). $$
Here $\Xx /Y$ is the category of pairs $(X,f)$ where $X\in \Xx$ and
$f: \pi (X)\rightarrow Y$ is a morphism in $\Yy$. Given a morphism
$a:Y\rightarrow Y'$ we get a functor $\alpha :\Xx /Y\rightarrow \Xx
/Y'$.   Furthermore
$$
F|_{\Xx /Y}= \alpha ^{\ast}(F|_{\Xx /Y'}),\;\;\;
G|_{\Xx /Y}= \alpha ^{\ast}(G|_{\Xx /Y'}).
$$
This and the remark of the previous paragraph gives a morphism
$$
\Re _{\Xx /Y}(F|_{\Xx /Y}, G|_{\Xx /Y})\rightarrow
\Re _{\Xx /Y'}(F|_{\Xx /Y'}, G|_{\Xx /Y'}),
$$
which is the morphism of functoriality for $\Re _{\Xx /\Yy}(F, G)$.

\begin{lemma}
\label{adjoint?}
Suppose $\pi : \Xx \rightarrow \Yy$ is a functor. Suppose that
$F: \Xx
\rightarrow Top$ is a covariant functor and $G: \Yy \rightarrow Top$ a
contravariant functor.  Then there is a natural weak homotopy
equivalence
$$
\Re _{\Xx} (F, \pi ^{\ast}(G))\cong
\Re _{\Yy}(\Re _{\Xx /\Yy}(F, \ast
),G).
$$
\end{lemma}
The proof is left to the reader since we shall not use it below.
\eop

Actually we don't use the above standard version; instead, we need
two more special versions of the relative realization. These are a
covariant version depending only on $F$ and a contravariant version
depending only on $G$.

Suppose that $\pi : \Xx \rightarrow \Yy$ is a split fibered category,
in other words it comes from a functor $\Phi : \Yy \rightarrow Cat$.
The objects of $\Xx$ are pairs of the form $(Y,U)$ with $U\in \Phi
(Y)$. Note that any $\Phi (Y)$ may be considered as a subcategory
of $\Xx$.
Suppose that $F: \Xx \rightarrow Top$ is a covariant
functor.  Then
we define the {\em special covariant relative realization} which is
a covariant functor denoted
$$
\Re _{\Xx /\Yy}^{\rightarrow}(F, \ast ): \Yy \rightarrow Top
$$
(note---in spite of the notation---that there is only one variable
$F$). It is defined by
$$
\Re _{\Xx /\Yy}^{\rightarrow}(F,\ast )(Y):=
\Re _{\Phi (Y)}(F|_{\Phi (Y)}, \ast _{\Phi (Y)}).
$$
The second variable $\ast_{\Phi (Y)}$ is
the constant functor on $\Phi (Y)$ associating to each object the
one point space.
If $a: Y\rightarrow Y'$ is a morphism in $\Yy$ then we get
a functor $\Phi (a): \Phi (Y)\rightarrow \Phi (Y')$.  The structure
of $F$ as functor on $\Xx$ gives a morphism of functors
$$
F|_{\Phi (Y)} \rightarrow \Phi (a)^{\ast}(F|_{\Phi (Y')}).
$$
We obtain the first of the following two morphisms:
$$
\Re _{\Phi (Y)}(F|_{\Phi (Y)}, \ast _{\Phi (Y)})
\rightarrow \Re _{\Phi (Y)}
(\Phi (a)^{\ast}(F|_{\Phi (Y')}), \ast _{\Phi (Y)})\rightarrow
\Re _{\Phi (Y')}(F|_{\Phi (Y')}, \ast _{\Phi (Y')}),
$$
the second being the standard morphism of functoriality
(noting obviously
that $\ast_{\Phi (Y)}= \Phi (a)^{\ast}(\ast_{\Phi
(Y')})$). The composition of the above two morphisms gives the
morphism of functoriality for $\Re _{\Xx /\Yy}^{\rm \rightarrow}(F,\ast )$.

\begin{lemma}
\label{covariant}
Suppose in the above situation that $G: \Yy \rightarrow Top$ is a
presheaf of spaces on $\Yy$. Then we have a natural equivalence
$$
\Re _{\Xx}(F, \pi ^{\ast}(G))\cong \Re _{Yy}(\Re ^{\rm
\rightarrow}_{\Xx /\Yy}(F, \ast ), G).
$$
\end{lemma}
{\em Proof:}
The realization on the right, with a realization in the argument, can be seen
as the realization of a bisimplicial space. If we take the diagonal simplicial
space this gives the realization on the left. The natural inclusion of the
realization of the diagonal into the realization of a bisimplicial space
is an equivalence.
\eop

We now get to the special contravariant version.  Suppose this time
that  $\pi : \Xx \rightarrow \Yy$ is a split fibered category
in the other direction
corresponding to a contravariant functor $\Phi : \Yy \rightarrow
Cat$.
The objects of $\Xx$ are again pairs of the form $(Y,U)$ with $U\in
\Phi (Y)$, and any $\Phi (Y)$ may be considered as a
subcategory of $\Xx$.
Suppose that $G: \Xx \rightarrow Top$ is a contravariant
functor.  Then
we define the {\em special contravariant relative realization} which
is a contravariant functor denoted
$$
\Re _{\Xx /\Yy}^{\leftarrow}(\ast ,G): \Yy \rightarrow Top
$$
(as before there is only one variable
$G$). It is defined by
$$
\Re _{\Xx /\Yy}^{\leftarrow}(\ast ,G)(Y):=
\Re _{\Phi (Y)}(\ast _{\Phi (Y)},G|_{\Phi (Y)}).
$$
The first variable $\ast_{\Phi (Y)}$ is
the constant functor on $\Phi (Y)$ associating to each object the
one point space.
If $a: Y\rightarrow Y'$ is a morphism in $\Yy$ then we get
a functor $\Phi (a): \Phi (Y')\rightarrow \Phi (Y)$.  The structure
of $G$ as presheaf on $\Xx$ gives a morphism of functors
$$
G|_{\Phi (Y')} \rightarrow \Phi (a)^{\ast}(G|_{\Phi (Y)}).
$$
We obtain the first of the following two morphisms:
$$
\Re _{\Phi (Y')}(\ast _{\Phi (Y')},G|_{\Phi (Y')})
\rightarrow
\Re _{\Phi (Y')}(\ast _{\Phi (Y')},\Phi (a)^{\ast}(G|_{\Phi (Y)})
\rightarrow
\Re _{\Phi (Y)}(\ast _{\Phi (Y)},G|_{\Phi (Y)}),
$$
the second being the standard morphism of functoriality
(noting as before
that $\ast_{\Phi (Y)}= \Phi (a)^{\ast}(\ast_{\Phi
(Y')})$). The composition of the above two morphisms gives the
morphism of functoriality for $\Re _{\Xx /\Yy}^{\leftarrow}(\ast
,G)$.

\begin{lemma}
\label{contravariant}
Suppose in the above situation that $F: \Yy \rightarrow Top$ is a
covariant functor. Then we have a natural equivalence
$$
\Re _{\Xx}(\pi ^{\ast}(F),G)\cong \Re _{Yy}(F,\Re ^{\rm
\leftarrow}_{\Xx /\Yy}(\ast ,G)).
$$
\end{lemma}
{\em Proof:}
This situation is a actually the same as the previous one after
taking the opposite categories, interchanging the roles of $F$ and
$G$ in the realization.
\eop

\numero{A descent condition for $F$}
In this and the remaining sections we assume that $\Xx$ has a
Grothendieck topology, i.e. $\Xx$ is a site.  Note however that the
definition of  $\Re _{\Xx}(F,G)$ does not depend on the Grothendieck
topology.

Suppose $F$ is a covariant functor $\Xx \rightarrow Top$. We say that
{\em $F$ satisfies covariant descent} if, for any object $X\in \Xx$
and any sieve $\Bb \subset \Xx /X$, the natural morphism
$$
\Re _{\Bb} (F|_{\Bb} , \ast ) \rightarrow \Re _{\Xx /X}(F|_{\Xx /X},
\ast )
$$
is a weak homotopy equivalence. Note that the space on the right is
equivalent to $F(X)$ by Lemma \ref{finalObject}.  The morphism
occuring above comes from our general discussion using the observation
that $F|_{\Bb}$ is the pullback from $\Xx /X$ to
$\Bb$ of $F|_{\Xx /X}$.

The following theorem gives the main example for our purposes.
\begin{theorem}
\label{example}
Suppose $\Xx$ is the site of schemes of finite type over $Spec (\cc )$
with the fppf topology (or any weaker topology such as the etale or
Zariski topology). The functor $F: \Xx \rightarrow Top$ defined by
setting $F(X)=X^{\rm top}$ (the usual topological space underlying the
analytic space associated to $X$) satisfies covariant descent.
\end{theorem}
{\em Proof:}
Suppose $\Bb$ is a sieve over $X$. Let $B$ denote the disjoint union of the
schemes in $\Bb$ with a morphism $B\rightarrow X$ (considered as a
non-noetherian scheme locally of finite type). This morphism is fppf surjective.
Let $N_{\cdot}(B/X)$ denote the standard simplicial scheme whose components are
$B\times _X \ldots \times _XB$.  Let $B^{\rm top}$ denote the associated
topological space (mapping to $X^{\rm top}$ and let $N_{\cdot}(B^{\rm
top}/X^{\rm top})$ again denote the simplicial space whose elements are fiber
products. There is no confusion in the notation because the functor $F$ commutes
with fiber products, so we can first take $B^{\rm top}$ then take the nerve, or
vice-versa getting the same answer.
We have
$$
\Re _{\Bb}(F|_{\Bb}, \ast )\cong
| N_{\cdot}(B^{\rm top}/X^{\rm top})|.
$$
It is not too hard to see (by stratifying everything and so on) that
$X^{\rm top}$ admits a triangulation as a simplicial complex which we denote
$S$ such that every simplex lifts into $B^{\rm top}$ (in the notation below
we will replace $X^{\rm top}$ by $S$).   Let $\tilde{S}$ denote the
simplicial complex which is the disjoint union of the simplices of $S$. We
can choose a map $\tilde{S}\rightarrow B^{\rm top}$.  On the other hand, the
morphism
$$
|N_{\cdot}(\tilde{S}/S)| \rightarrow S
$$
is a Serre fibration and a
weak equivalence.  Consider the bisimplicial space
$N_{i}(B^{\rm top}/S) \times _S N_{j }(\tilde{S}/S)$.
If we  fix the variable $j$
then in the variable $i$ it is the same as the nerve of the map
$$
B^{\rm top} \times _S N_j (\tilde{S}/S)\rightarrow N_j(\tilde{S}/S).
$$
Since this map admits a section, the realization of its nerve is weakly
equivalent to the base. Thus if we realize first in the $i$ direction and then
in the $j$ direction we obtain something mapping by a weak equivalence to $S$.
On the other hand if we fix $i$ then the realization in the $j$ variable is just
$$
N_i(B^{\rm top} /S)\times _S|N_{\cdot}(\tilde{S}/S),
$$
which maps by a weak equivalence to $N_i(B^{\rm top}/S)|$.  Thus if we realize
first in the $j$-direction and then in the $i$-direction we obtain something
mapping by a weak equivalence to
$| N_{\cdot}(B^{\rm top}/X^{\rm top})|$. This proves that the map
$$
| N_{\cdot}(B^{\rm top}/X^{\rm top})|\rightarrow X^{\rm top}
$$
is a weak equivalence.
\eop

We have a similar analytic version.  Let $\Xx ^{\rm an}$ denote the
site whose underlying category is that of complex analytic spaces, and
whose topology is given by saying that a family is surjective if it
admits sections locally on the base.  Let $F: \Xx ^{\rm an}\rightarrow
Top$ be the functor associating to an analytic space, its underlying
topological space.

\begin{theorem}
\label{analyticExample}
With the above notations, $F: \Xx ^{\rm an}\rightarrow Top$ satisfies
covariant descent.
\end{theorem}
{\em Proof:}
Use the same proof as above but with $S$ just as $X^{\rm top}$ and $\tilde{S}$
as an open covering of $S$ admitting a lifting to $B^{\rm top}$.
The main step, that $|N_{\cdot}(\tilde{S}/S)|\rightarrow S$ is a Serre
fibration and weak equivalence, still  holds.
\eop

{\em Remark:} If we just want Theorem \ref{example} for the etale or Zariski
topologies then we can use the ``open'' version as in \ref{analyticExample},
this avoids the use of a triangulation.

\numero{Statement of the main theorem}

Suppose as above that $\Xx$ is a site. Recall that a morphism
$G\rightarrow G'$ of presheaves of topological spaces (i.e.
contravariant functors $\Xx \rightarrow Top$) is called an {\em
Illusie weak equivalence} \cite{Illusie} if it induces isomorphisms
of the {\em sheaves associated to the homotopy presheaves $\pi
_0(G)$ or $\pi _i (G, g)$} (this is explained further in the
introduction).

\begin{theorem}
\label{main}
Suppose $\Xx$ is a site and $F: \Xx \rightarrow Top$ is a covariant
functor which satisfies covariant descent.  Then any Illusie weak
equivalence $G\rightarrow G'$ between contravariant functors induces a
weak homotopy equivalence of realizations
$$
\Re _{\Xx}(F,G)\stackrel{\cong}{\rightarrow} \Re _{\Xx }(F, G').
$$
\end{theorem}

Even though the theorem can be stated without refering to the
closed model category structure or homotopy sheafification, these
notions are essential in our proof.  Actually I think that there
also exists a proof which doesn't use these ideas but which
proceeds by Postnikov  induction. However, one needs to start at
degree zero with the same result for presheaves of sets, and as it
seems that the argument needed for presheaves of sets is essentially
the same as our argument below, we have prefered to go straight
through with this argument in general.
The following section treats the generalities we will need, the
subsequent section gives the main lemma.

\numero{Homotopy-sheafification}
Suppose $G$ is a presheaf of spaces on a category $\Cc$ (i.e. a
contravariant functor to $Top$). A {\em section of $G$ over $\Cc$}
is a function $g$ which to each composable sequence
$\phi _{\cdot} : X_0 \rightarrow X_n$
associates $g (\phi _{\cdot}): \Delta ^n \rightarrow G(X_0)$
such that if $\eta : \Delta ^m\rightarrow \Delta ^n$ is a face,
and if $\phi '_{\cdot} : X_i\rightarrow X_j$ is the composable
sequence obtained by applying the face to $\phi_{\cdot}$, then
$$
g(\phi _{\cdot} )\circ \eta = g(\phi '_{\cdot}).
$$
The space of sections  can be topologized using the compact-open
topology on the space of maps $\Delta ^n\rightarrow G(X_0)$, taking
the Tychonoff topology on the product over all $\phi _{\cdot}$ of
these spaces of maps, and then considering the space of sections as a
subspace of the product.  Let $\Gamma (\Cc , G)$ denote this space of
sections.

If $\pi : \Cc \rightarrow \Dd$ is a functor and if $G$ is a presheaf
of spaces over $\Dd$ then there is an induced map $\Gamma (\Dd ,
G)\rightarrow \Gamma (\Cc , \pi ^{\ast}(G))$.

If $\Cc$ has a final object $U$ then by Lemma \ref{finalObject}
the
map $\Gamma (\Cc , G)\rightarrow G(U)$ is a weak homotopy
equivalence.

Now we get back to our site $\Xx$.  For any object $X\in \Xx$ we have
the weak equivalence $\Gamma (\Xx /X, G|_{\Xx /X})\cong G(X)$. If 
$\Bb \subset \Xx /X$ is a sieve then we get a natural map
$$
\Gamma (\Xx /X, G|_{\Xx /X})\rightarrow \Gamma (\Bb , G|_{\Bb}).
$$
We say that $G$ is a {\em homotopy-sheaf} if for any object $X\in \Xx$
and any sieve $\Bb\subset \Xx /X$ this restriction map is a weak
homotopy equivalence.

Recall from \cite{Jardine1} that for any presheaf of spaces $G$ (which
can be considered as a simplicial presheaf by taking the singular
simplicial set over each object, for example) there is an essentially
canonical morphism $G\rightarrow G'$ to a {\em fibrant} presheaf of
spaces. While not getting into the precise definition of Jardine's
fibrant condition here, we note that the essential point is that a
fibrant simplicial presheaf corresponds to a presheaf of spaces which
is a homotopy sheaf. To prove this we can interpret $\Gamma (\Bb ,
G|_{\Bb})$ as a morphism space $Hom (\ast '_{\Bb}, G)$ where
$\ast '_{\Bb}$ is a presheaf equivalent to the presheaf of sets
$\ast_{\Bb}$ represented by $\Bb$, i.e. it gives a contractible space
for each object of $\Bb$ and empty otherwise.  The presheaf of spaces
$\ast '_{\Bb}$ is Illusie weak-equivalent to $\ast _{\Xx /X}$ so (by
the closed model structure \cite{Jardine1}) the morphism
$$
Hom (\ast '_{\Bb}, G)\rightarrow Hom (\ast _{\Xx /X}, G)
$$
is a weak homotopy equivalence of spaces.

\begin{proposition}
\label{uniqueness}
The replacement  of $G$ by a homotopy-sheaf $G'$ in the $n$-truncated case
(resp.
fibrant object in the non-truncated case) together with an  Illusie weak
equivalence $G\rightarrow G'$ is unique up to object-by-object equivalence. More
precisely if $G' $ and $ G''$ are two such replacements
with Illusie weak equivalences $G\rightarrow G'$ and $G\rightarrow G''$
then
there is a diagram
$$
\begin{array}{ccc}
G & \rightarrow G'' \\
\downarrow && \downarrow \\
G' & \rightarrow & G^{(3)}
\end{array}
$$
where the bottom horizontal arrow and the right vertical arrow are
object-by-object equivalences.  Furthermore this diagram is unique up
to homotopy.
\end{proposition}
{\em Proof:}
We can assume that $G'$ and $G''$ are fibrant in the sense of Jardine,
for in the $n$-truncated case if $H$ is a homotopy-sheaf then the replacement
by a fibrant object $H\rightarrow H'$ is an object-by-object weak equivalence.

(To see this note that it suffices to show that $H\times _{H'}Z\rightarrow Z$
is an object-by-object equivalence for any $Z\in \Xx$ and $Z\rightarrow H'$;
but this fiber product preserves the homotopy-sheaf condition so we may
effectively assume that $H'$ is represented by $Z$; then arguing inductively on
$n$ we show that the top homotopy-group presheaf is trivial so $H$ is
$n-1$-truncated, eventually we get down to $H$ $0$-truncated and then an
Illusie weak equivalence between sheaves of sets is an isomorphism.)

Let
$F$ be the pushout of $G\rightarrow G'$ and $G\rightarrow G''$, then let
$F\rightarrow G^{(3)}$ be a replacement by a fibrant object.  The morphism
$G'\rightarrow G^{(3)}$ is an Illusie weak equivalence between
fibrant objects.  Therefore this morphism is invertible up to homotopy,
in particular it is an object-by-object equivalence. The same goes for $G''
\rightarrow G^{(3)}$.
\eop

If $G'$ is a homotopy sheaf and $G\rightarrow G'$ is an Illusie weak
equivalence then we say that $G'$ is the {\em homotopy-sheafification
of $G$}.

Suppose now that $G$ is {\em $n$-truncated} (i.e. the homotopy groups
of the $G(X)$ vanish in degrees strictly bigger than $n$).  Then we
propose an explicit method for obtaining the homotopy sheafification
by applying $n+2$ times a certain process denoted $\gamma$.  In the
case of presheaves of sets ($n=0$) this is essentially the same thing
as the standard construction of the sheafification of a presheaf of
sets by doing the obvious operation $2$ times.   We will define
$\gamma (G)$ below with a morphism $G\rightarrow \gamma (G)$ (an
Illusie weak equivalence).  It was shown in \cite{flexible} that if
$G$ is $n$-truncated then $\gamma ^{n+2}G$ is a homotopy sheaf. We
give a sketch of proof below also, since \cite{flexible} is not
published.

Note that Jardine's proof of the closed
model structure \cite{Jardine1} gives, of course, a method for
obtaining $G\rightarrow G'$ with $G'$ fibrant (and thus obtaining the
homotopy-sheafification). His method, while explicit, is not very
useful for proving things about the notion of Illusie weak equivalence
because it consists of taking the pushout over {\em all} diagrams
$U\stackrel{a}{\leftarrow} V\rightarrow G$ with $a$ an Illusie weak
equivalence (and then iterating this an infinite number of times using
transfinite induction).

Using our operation $\gamma$ we have the following outline of the
proof of Theorem \ref{main}.  The main lemma is to prove, using
the explicit definition of $\gamma$, that the morphism
$G\rightarrow \gamma (G)$ induces a weak equivalence of
realizations
$$
\Re _{\Xx}(F,G)\rightarrow \Re _{\Xx}(F , \gamma (G)).
$$
The same is of course true for the iteration $G\rightarrow \gamma
^{n+2}(G)$.
Now suppose $G\rightarrow G'$ is an Illusie weak equivalence
between $n$-truncated presheaves of spaces.  Then we obtain a diagram
$$
\begin{array}{ccc}
G & \rightarrow & \gamma ^{n+2}(G) \\
\downarrow && \downarrow \\
G' & \rightarrow & \gamma ^{n+2}(G')
\end{array}
$$
giving a corresponding diagram on the level of realizations.
The horizontal arrows induce weak equivalences of the realizations.
On the other hand, all arrows are Illusie weak equivalences and the
spaces on the right are homotopy-sheaves. Thus by the uniqueness result,
the right vertical arrow is an object-by-object equivalence, in
particular it induces a weak equivalence of realizations.  We conclude
that the vertical arrow on the left induces a weak equivalence of
realizations, which is Theorem \ref{main}.

To complete the proof of Theorem \ref{main} we will just need a
little argument to go from the $n$-truncated case to the general
case (since the hypothesis of $n$-truncation is not present in the
statement of \ref{main}).  We give this below.

\subnumero{The definition of $\gamma$}

Suppose $\Xx$ is a site and $G$ is a presheaf of spaces on $\Xx$.
Let $\pi : \Ee \rightarrow \Xx$ be the fibered category associated
to the contravariant functor $X\mapsto Sv(X)$ where $Sv(X)$ denotes
the category (directed set, really) of sieves $\Bb \subset \Xx /X$.
The objects of $\Ee$ are pairs $(X,\Bb )$ where $X\in \Xx$ and
$\Bb \subset \Xx /X$ is a sieve.

We define the contravariant functor $\tilde{\gamma} (G): \Ee
\rightarrow Top$ by
$$
\tilde{\gamma}(G)(X, \Bb ):= \Gamma (\Bb , G|_{\Bb}).
$$
We are in the situation of the contravariant relative realization
discussed at the start; we  put
$$
\gamma '(G):= \Re ^{\leftarrow}_{\Ee /\Xx}(\ast , \gamma '(G)).
$$
On the other hand note that there is a canonical section
$$
\xi: \Xx \hookrightarrow \Ee
$$
defined by setting $\xi (X):= (X, \Xx /X)$. Put
$$
\sigma (G):= \xi ^{\ast}(\tilde{\gamma}(G)).
$$
Thus
$$
\sigma (G)(X):= \Gamma (\Xx /X, G).
$$
There is a natural weak equivalence
$$
\sigma (G)\stackrel{\rightarrow} G.
$$
On the other hand note trivially
that
$$
\sigma (G)=\Re ^{\leftarrow} _{\xi (\Xx )/\Xx }(\ast ,
\tilde{\gamma}(G) |_{\xi (\Xx )}),
$$
and from this and a functoriality of the realization in terms of
the base category we obtain a natural morphism
$$
\sigma (G)\rightarrow \gamma ' (G).
$$
Finally let $\gamma (G)$ denote the push-out (object-by-object)
of the diagram
$$
\begin{array}{ccc}
\sigma (G) & \rightarrow & \gamma ' (G)\\
\downarrow && \\
G && .
\end{array}
$$
We obtain a morphism $\gamma (G)$ which is object-by-object weak
equivalent to the morphism $\sigma (G)\rightarrow \gamma '(G)$.

{\em Remark:} We have
$$
\gamma '(G)(X) = \lim _{\rightarrow , \Bb \subset \Xx /X}
\Gamma (\Bb , G|_{\Bb}).
$$
In particular, $\gamma '(G)$ (or, up to object-by-object weak
equivalence, $\gamma (G)$) is the same as the presheaf of spaces
constructed in \cite{flexible}.

{\em Remark:} The operation $\gamma$ is functorial in $G$. This is
completely clear from the definition since at each step we apply a
functor.

{\em Remark:} If $G$ is $n$-truncated then so is $\gamma (G)$.

\subnumero{The homotopy-sheafification $\gamma ^{n+2}$}

\begin{theorem}
Suppose $G$ is an $n$-truncated presheaf of spaces. Then
the morphism obtained by iterating the operation $\gamma$
$$
G\rightarrow \gamma ^{n+2}(G)
$$
is the homotopy-sheafification of $G$.  In Jardine's terms this
means that the morphism is object-by-object weak equivalent to the
replacement of $G$ by a fibrant object.
\end{theorem}
{\em Proof:}
For now we refer to \cite{flexible} for the proof.
\eop

\numero{The main lemma}

We treat the situation we will encounter in the main lemma, in a
general context first.

Suppose $\Cc$ is a category with a projector,
that is a functor $P: \Cc \rightarrow \Cc$ such that $P^2= I$ (where
$I$ denotes the identity functor of $\Cc$).  Let $\Dd = P(\Cc )
\subset \Cc$ denote the image subcategory. Suppose we have a natural
morphism of functors $\psi : P\rightarrow I$ such that $\psi
|_{\Dd}$ is the identity of $P|_{\Dd}= I_{\Dd}$. Suppose that $F:
\Dd \rightarrow Top$ is a covariant functor and $G: \Cc \rightarrow
Top$ is a contravariant functor (presheaf of spaces).  Then $\psi $
induces a morphism of presheaves of spaces $G\rightarrow
P^{\ast}(G)$. On the other hand,  $P$ induces a morphism of
realizations $\Re _{\Cc}(P^{\ast}(F),
P^{\ast}(G|_{\Dd}))\rightarrow  \Re _{\Dd}(F,G|_{\Dd})$.   Composing
these we obtain a morphism $$ a: \Re _{\Cc}(P^{\ast}(F),
G)\rightarrow \Re _{\Dd}(F, G|_{\Dd}). $$
On the other hand, note that $P^{\ast}(F)|_{\Dd} = F$. The inclusion
$\Dd \subset \Cc$ induces a morphism
$$
b: \Re _{\Dd} (F, G|_{\Dd}) \rightarrow \Re _{\Cc}(P^{\ast}(F), G).
$$
\begin{lemma}
\label{projector}
With the above notations, $ab$ is equal to the identity and $ba$ is
homotopic to the identity.  In particular $a$ and $b$ are homotopy
equivalences of spaces.
\end{lemma}
{\em Proof:}
The hard part is the case of $ba$. Note that
$$
\Re _{\Cc}(P^{\ast}(F), G)= \coprod _{\phi _{\cdot} : X_0\rightarrow
X_n \in \Cc} F(P(X_0)) \times G(X_n) \times \Delta ^n /\sim .
$$
Given a point $(\phi _{\cdot}, f,g,t)$ here, the image by  $a$ is the
point
$(P(\phi _{\cdot}), f,\psi _{X_n}^{\ast}(g),t)$ in the component
$F(P(X_0)) \times G(P(X_n)) \times \Delta ^n$ of
$$
\Re _{\Dd}(F, G|_{\Dd})= \coprod _{\rho _{\cdot} :
Y_0\rightarrow Y_n \in \Dd} F(Y_0) \times G(Y_n) \times \Delta ^n
/\sim .
$$
Since $b$ is just induced by the inclusion, the image of this point by
$b$ can be written in the same way as
$(P(\phi _{\cdot}), f,\psi _{X_n}^{\ast}(g),t)$, a point in the
component
$F(P(X_0)) \times G(P(X_n)) \times \Delta ^n$ of
$$
\Re _{\Cc}(P^{\ast}(F), G)= \coprod _{\rho _{\cdot} :
Y_0\rightarrow Y_n \in \Cc} F(Y_0) \times G(Y_n) \times \Delta ^n
/\sim .
$$
To conclude, we can write
$$
ba(\phi _{\cdot}, f, g, t)=
(P(\phi _{\cdot}), f,\psi _{X_n}^{\ast}(g),t).
$$
We would like to define a homotopy from $ba$ to the identity.  To do
this note that the realization of the diagram
$$
\begin{array}{ccccccc}
P(X_0) & \rightarrow & P(X_1)& \rightarrow & \ldots & \rightarrow &
P(X_n) \\
\downarrow {\displaystyle \psi _{X_0}}& & \downarrow
{\displaystyle \psi _{X_1}} &    &                &   & \downarrow
{\displaystyle \psi _{X_n}}\\ X_0 & \rightarrow & X_1& \rightarrow &
\ldots & \rightarrow & X_n \end{array}
$$
is naturally isomorphic to $[0, 1]\times \Delta ^n$ (this is an
example of the division of the product of simplices which can be
accomplished via the nerves of categories).  From this diagram  we
obtain a map
$$
H_{\phi _{\cdot}}: [0,1]\times \Delta ^n \times F(P(X_0))\times G(X_n)
\rightarrow \Re _{\Cc}  (P^{\ast}(F), G),
$$
such that
$$
H_{\phi _{\cdot}}(0, t,f,g)= (P(\phi _{\cdot}), f, \psi
_{X_n}^{\ast}(g), t)
$$
and
$$
H_{\phi _{\cdot}}(1, t,f,g)= (\phi _{\cdot}, f, g,t).
$$
We leave it to the reader to check that these maps fit together with
the glueing relations $\sim$ to give a map
$$
H: [0,1]\times \Re _{\Cc}(P^{\ast}(F),G)\rightarrow
\Re _{\Cc}(P^{\ast}(F),G)
$$
such that $H(0) = ba$ and $H(1)$ is the identity.

It is easy to see that $ab$ is the identity, using the above
formula for $a$ and the fact that $\psi |_{\Dd}$ is the identity on
$P|_{\Dd}=I_{\Dd}$.
\eop

\medskip

We can now state and prove the main lemma.

\begin{lemma}
\label{mainLemma}
Suppose $G$ is a presheaf of spaces, and suppose $F$ is a covariant
functor $\Xx \rightarrow Top$ satisfying covariant descent.  Then the
morphism $G\rightarrow \gamma (G)$ induces a weak equivalence of
realizations $\Re _{\Xx}(F,G)\cong \Re _{\Xx }(F,\gamma (G))$.
\end{lemma}
{\em Proof:}
Let $\Cc$ be the category of triples $(X, \Bb , Y)$ where $X\in \Xx$,
$\Bb \subset \Xx /X$ is a sieve, and $Y\in \Bb$ is an object (together
with morphism to $X$).  We have an inclusion of categories $i: \Xx
\hookrightarrow \Cc$ defined by $i(X)= (X, \Xx /X, X)$.
Let $\Dd$ be the image of $i$ (it is equal to $\Xx$).  Let $P: \Cc
\rightarrow \Cc$ be the morphism defined by
$$
P(X, \Bb , Y):= i(Y).
$$
The morphism $Y\rightarrow X$ (which comes with the data of $Y$)
gives a morphism $i(Y)\rightarrow (X,\Bb , Y)$. This is natural, so it
provides our natural transformation $\psi : P\\rightarrow I$.
Define a presheaf of spaces $G'$  on $\Cc$ by
$$
G'(X, \Bb , Y):= \Gamma
(\Bb , G|_{\Bb}).
$$
Define the presheaf of spaces $F'$ on $\Cc$ by
$$
F'(X, \Bb , Y):= F(Y).
$$
Note that $F' = P^{\ast}(F')$ already.
Apply Lemma \ref{projector} to $F'$ and $G'$ on $\Cc$ with $P$ and
$\psi$ as above.  We obtain that the inclusion
induces a homotopy equivalence (denoted $b$ in Lemma \ref{projector})
$$
\Re _{\Xx}(i^{\ast}(F'),
i^{\ast}(G'))
=
\Re _{\Dd}(F',
G'|_{\Dd})
\rightarrow
\Re _{\Cc}(P^{\ast}(F'), G')\stackrel{\cong}{\rightarrow}
\Re _{\Dd}(F',
G'|_{\Dd}).
$$
We now interpret these elements in terms of $F$ and $G$.
The left hand side first.  Note that $i^{\ast}(F')=F$ whereas
$i^{\ast}(G')= \sigma (G)$ in the notations used in the definition
of $\gamma$.
We obtain
that the inclusion $i$ induces a weak homotopy equivalence
$$
\Re _{\Xx}(F, \sigma (G))
\stackrel{\cong}{\rightarrow}
\Re _{\Cc}(P^{\ast}(F'), G') .
$$
Next, define an intermediate category $\Ee$ to be the category of
pairs $(X,\Bb )$, with a functor $q: \Ee \rightarrow \Xx$ defined by
$q(X,\Bb )= X$.  This is the same as the category used to define
$\gamma '$ and $\gamma$. The presheaf $G'$ is pulled back from
the one we
denoted $\tilde{\gamma}(G)$ on $\Ee$.  We have
$$
\Re
_{\Cc}(P^{\ast}(F'), G')= \Re _{\Ee}(\Re ^{\rightarrow}_{\Cc
/\Ee}(P^{\ast}(F'), \ast ), \tilde{\gamma}(G)).
$$
There is
a natural morphism
$$
\mu : \Re ^{\rightarrow}_{\Cc /\Ee}(P^{\ast}(F'),\ast )\rightarrow
q^{\ast}(F).
$$
The condition that $F$ satisfies
covariant descent is exactly the condition that $\mu$ is an
object-by-object weak equivalence.
We obtain a natural weak equivalence
$$
c:\Re _{\Cc}(P^{\ast}(F'), G') \stackrel{\cong}{\rightarrow}
\Re _{\Ee}(q^{\ast}(F), \tilde{\gamma}(G) ).
$$
The map $i$ composed with the projection to $\Ee$ is the map
$\xi : \Xx \rightarrow \Ee$ defined by $\xi (X)= (X, \Xx /X)$.  Note
that $\xi ^{\ast}(q^{\ast}(F)) = F$ whereas
$\xi ^{\ast}(\tilde{\gamma}(G))=
\sigma (G)$.
Also $\xi ^{\ast}(P^{\ast}(F))= F$ and the composition
$$
F\rightarrow \xi ^{\ast}\Re ^{\rightarrow}_{\Cc
/\Ee}(P^{\ast}(F'),\ast ) \rightarrow F
$$
is the identity, so the composition of our above morphism
$c$ with the
inclusion $b$ gives the map
$$
\Re _{\Xx}(F, \sigma (G))\rightarrow \Re _{\Ee}(q^{\ast}(F),
\tilde{\gamma}(G) )
$$
induced by $\xi $. Again this map is a weak equivalence (the fact
that $b$ is a weak equivalence is Lemma \ref{projector} while the
fact that $c$ is a weak equivalence is the covariant descent
condition for $F$).

Finally we have a natural weak equivalence
$$
\Re _{\Ee}(q^{\ast}(F),
\tilde{\gamma}(G) )\cong \Re _{\Xx}(F, \Re ^{\leftarrow}_{\Ee /\Xx}
(\ast , \tilde{\gamma}(G) )).
$$
By definition
$$
\Re ^{\leftarrow}_{\Ee /\Xx}
(\ast , \tilde{\gamma}(G) ) = \gamma '(G),
$$
and the morphism $\sigma (G)\rightarrow \Re ^{\leftarrow}_{\Ee /\Xx}
(\ast , \tilde{\gamma}(G) )$
induced by $\xi $ is the morphism constructed above.
We conclude by composing this with the previous morphism to 
obtain the
weak equivalence
$$
\Re _{\Xx}(F, \sigma (G))\rightarrow \Re _{\Xx}(F,\gamma '(G)).
$$
Finally note that the pushout diagram for $\gamma (G)$ gives
(by functoriality of the realization) a diagram
$$
\begin{array}{ccc}
\Re _{\Xx}(F, \sigma (G))&\rightarrow &\Re _{\Xx}(F,\gamma '(G))\\
\downarrow & & \downarrow \\
\Re _{\Xx}(F, G)&\rightarrow &\Re _{\Xx}(F,\gamma (G)).
\end{array}
$$
The vertical arrows are obtained from object-by-object weak
equivalences in the pushout diagram of $\gamma (G)$, so they are weak
equivalences of spaces here.  Thus the bottom arrow is a weak
equivalence.
This completes the proof of
the main lemma. \eop

\subnumero{Proof of Theorem \ref{main}}

We have treated above the case where $G$ is $n$-truncated.  For any $G$ we have
the morphism $G\rightarrow \tau _{\leq n}^{\rm pre}G$ truncating the homotopy
groups below $n$, object-by-object.  if $G\rightarrow G'$ is an Illusie weak
equivalence then for any $n$,
$\tau _{\leq n}^{\rm pre}G\rightarrow \tau _{\leq n}^{\rm pre}G'$ is an
Illusie weak equivalence of $n$-truncated presheaves of spaces. By the main
theorem for the truncated case,
$$
\Re _{\Xx}(F, \tau _{\leq n}^{\rm pre}G)\rightarrow
\Re _{\Xx}(F, \tau _{\leq n}^{\rm pre}G ')
$$
is a weak equivalence.  Note now that the morphism $G\rightarrow
\tau _{\leq n}^{\rm pre}G$ may be obtained by adding on balls of dimension
$\geq n+2$.  Hence the morphism
$$
\Re _{\Xx}(F, G)\rightarrow
\Re _{\Xx}(F, \tau _{\leq n}^{\rm pre}G )
$$
is obtained by adding on pieces which are products of balls of dimension $\geq
n+2$ with something else. Thus for $i\leq n$ the above morphism induces
isomorphisms
$$
\pi _i(\Re _{\Xx}(F, G) , y)\stackrel{\cong}{\rightarrow}
\pi _i(\Re _{\Xx}(F, \tau _{\leq n}^{\rm pre}G ),y).
$$
Combining with the result of the theorem for the truncated case we get that
$$
\pi _i(\Re _{\Xx}(F, G) , y)\stackrel{\cong}{\rightarrow}
\pi _i(\Re _{\Xx}(F, G') , y')
$$
for any $i\leq n$. Since this is true now for any $n$, we get that
$\Re _{\Xx}(F, G)\rightarrow \Re _{\Xx}(F, G')$ is a weak equivalence, which is
the conclusion of Theorem \ref{main}.
\eop

\numero{Applications}

\begin{lemma}
\label{disjoint}
Suppose $F: \Xx\rightarrow Top$ is a covariant functor and suppose
$G:\Xx \rightarrow Sets \subset Top$ is a presheaf of sets represented
(as a presheaf) by a disjoint union of objects of $\Xx$,
$$
G= \coprod
_{i\in I} X_i.
$$
Then
$$
\Re _{\Xx }(F, G) \cong \coprod _{i\in I} F(X_i).
$$
\end{lemma}
{\em Proof:}
It is immediate that
$$
\Re _{\Xx }(F, G) \cong \coprod _{i\in I} \Re _{\Xx }(F, X_i).
$$
On the other hand, by Lemma \ref{finalObject} we have that
$\Re _{\Xx }(F, X_i)\cong X_i ^{\rm top}$.
\eop

\begin{lemma}
\label{simplicial}
Suppose $F: \Xx \rightarrow Top$ is a covariant functor and suppose
$A_{\cdot}$ is a simplicial presheaf of sets on $\Xx$.  Let  $\Re
_{\Xx}(F, A_{\cdot})$ be the simplicial space whose components are the
realizations of the component presheaves of $A_{\cdot}$.  Let $G$
be the realization of $A_{\cdot}$ (i.e. $G(X)$ is the realization
of the simplicial set $A_{\cdot}(X)$ for each object $X$).
Then
$$
\Re _{\Xx}(F,G) = \Re (\Re _{\Xx}(F, A_{\cdot})).
$$
\end{lemma}
{\em Proof:}
Both sides are the realizations of the same bisimplicial set.
\eop

The preceding two lemmas allow us to identify $\Re _{\Xx}(F, G)$ with
the naive construction one would make: we can express $G$ (up to
object-by-object equivalence) as the realization of a simplicial
presheaf of sets, and we can even suppose that the component
presheaves are disjoint unions of objects of $\Xx$; then we can take
the realization of each component as simply the disjoint union of the
$F(X_i)$ with $X_i$ making up the component; and finally we can take
the realization of the resulting simplicial topological space.  The
two previous lemmas show that this process actually gives $\Re
_{\Xx}(F,G)$ up to weak homotopy equivalence.  Our main result now
shows (under the hypothesis that $F$ satisfies covariant descent) that
this construction is unchanged if we replace $G$ by an
Illusie-equivalent $G'$.  This doesn't seem to be easy to see directly
from the point of view of the construction given in this paragraph.

In view of the discussion above, we propose an alternative shorter
notation.  If the site $\Xx$ in question is understood from context,
then we denote by $F(G)$ the realization $\Re _{\Xx}(F,G)$. This is
compatible with the original notation for $F$ by Lemma
\ref{disjoint}.

We can apply this discussion to $1$-stacks.  An algebraic stack $\Gg$
is typically given by a {\em presentation}, starting with a surjective
morphism $X\rightarrow \Gg$ from a scheme of finite type $X$ (this
being a surjection of stacks but not pre-stacks!). The definition of
algebraic stack is such that $R:= X\times _{\Gg}X$ is represented by a
scheme of finite type.  We obtain a simplicial scheme
$$
Z_n := X\times _{\Gg} \ldots \times _{\Gg} X,
$$
with $Z_0=X$ and $Z_1=R$.
In the case of stacks over $Spec (\cc )$ we would like to define the
{\em realization} $|\Gg |$ to be the realization of the simplicial
topological space $\{ Z_n^{\rm top} \}$.  The realization of the
simplicial presheaf $|Z_{\cdot}|$ is a presheaf of spaces which is
Illusie weak equivalent to the $1$-truncated presheaf of spaces $G$
corresponding to the stack $\Gg$ (in fact, the presheaf of spaces
corresponding to $\Gg$ is the fibrant object associated to
$|Z_{\cdot}|$).  Our main theorem coupled with Lemmas \ref{disjoint}
and \ref{simplicial} shows that the realization of the simplicial
space  $\{ Z_n^{\rm top} \}$ is weak homotopy equivalent to the
realization $\Re _{\Xx }(F,G )$ as defined above (with $F$ as in
Theorem \ref{example}). In particular it is independant of the choice
of surjection $X\rightarrow \Gg$.

\numero{Algebraic to analytic}

Let $\Xx$ denote the category of schemes of finite type over $Spec
(\cc )$ and let $\Xx ^{\rm an}$ denote the category of
complex analytic spaces.  Suppose $G$ is a presheaf of spaces on
$\Xx$. We will define a presheaf of spaces $G^{\rm an}$ on $\Xx
^{\rm an}$, extending the construction of the complex analytic space
associated to a scheme of finite type.

For any $Z\in \Xx ^{\rm an}$ put
$$
G^{\rm an, pre}(Z):= \lim _{\rightarrow , B\subset \Gamma (Z, \Oo )}
G(Spec(B)).
$$
The limit is a homotopy limit (which can be thought of as a
realization over the indexing category, for example).  The limit is
taken over all subalgebras $B\subset \Gamma (Z, \Oo )$ of finite type
over $\cc$.
Let $G^{\rm an}$ be the homotopy sheaf associated to $G^{\rm an, pre}$.

Suppose $G$ is a sheaf of sets represented by a disjoint union
of schemes of finite type
$$
G= \coprod _{i\in I}X_i .
$$
Then
$$
G^{\rm an, pre} \cong  \coprod _{i\in I}X^{\rm an}_i
$$
that is $G^{\rm an}$ is (homotopic to) the sheaf of sets
represented by  the disjoint union of the associated analytic spaces.

If $\{ A_{\cdot}^i\} _{i\in I}$ is a directed system of simplicial
presheaves of spaces, then
$$
\lim _{\rightarrow , i\in I} |A_{\cdot}^i| \cong |(\lim _{i\in
I}A_{\cdot}^i)|,
$$
where $|A_{\cdot}|$ denotes the object-by-object realization of a
simplicial presheaf of spaces into a presheaf of spaces. If
$B_{\cdot}$ is a simplicial presheaf of spaces and if $B'_{\cdot}$
denotes the homotopy sheafification at each stage, then
$$
|B_{\cdot}|\rightarrow |B'_{\cdot }|
$$
is an Illusie weak equivalence ({\em proof:} it is the realization of
a morphism of simplicial presheaves of spaces which is an Illusie
weak equivalence at each stage; by the cohomological interpretation
of Illusie weak equivalence \cite{Illusie} we obtain the stated
result) and in particular induces an equivalence of homotopy
sheafifications.  Putting these together we find that if $G_{\cdot}$
is a simplicial presheaf of spaces on $\Xx$ then $(|G_{\cdot}|
)^{\rm an}$ is the homotopy sheafification of $|G^{\rm
an}_{\cdot}|$. We can apply this in the case where $A_{\cdot}$ is a
simplicial presheaf of sets whose components $G_n$ are represented
by disjoint unions of schemes of finite type, say $$ A_n = \coprod
_{i\in I_n}X_{n, i}. $$ Let $G= |A_{\cdot}|$ be the realization into
a presheaf of spaces. We find that $G^{\rm an}$ is the homotopy
sheafification of $A^{\rm an}_{\cdot}$ with
$$
A^{\rm an}_n = \coprod _{i\in I_n}X^{\rm an}_{n, i}
$$

We can apply the preceding discussion to obtain a comparison theorem
between the realizations of $G$ and $G^{\rm an}$.

\begin{theorem}
\label{comparison}
Suppose $F: \Xx \rightarrow Top$ and $F: \Xx ^{\rm an}\rightarrow
Top$ are the covariant functors given in Theorems \ref{example} and
\ref{analyticExample} (no confusion results from keeping the same
notation for the two).
Then for any presheaf of spaces $G$ on $\Xx$ we
have a  weak homotopy equivalence
$$
\Re _{\Xx} (F, G) \cong \Re _{\Xx ^{\rm an}}(F, G^{\rm an}).
$$
\end{theorem}
{\em Proof:}
Any presheaf of spaces $G$ is Illusie weak  equivalent (or in fact
object-by-object weak equivalent) to the realization of a presheaf of
sets $A_{\cdot}$ of the form discussed above (the components being
disjoint unions of schemes). On the other hand, we have
$$
\Re
_{\Xx}(F, |A_{\cdot}|)\cong  |\Re _{\Xx}(F, A_{\cdot})|, $$
and similarly
$$
\Re _{\Xx^{\rm an}}(F, |A^{\rm an}_{\cdot}|)\cong
|\Re _{\Xx}(F, A^{\rm an}_{\cdot})|.
$$
Finally, if $A_n= \coprod _{i\in I_n}X_{i,n}$ then
$$
\Re _{\Xx}(F, A_n )= \coprod _{i\in I_n}X_{i,n}^{\rm top} =
\Re _{\Xx^{\rm an}}(F, A^{\rm an}_n ).
$$
We obtain that
$$
\Re _{\Xx}(F, |A_{\cdot}|) \cong
\Re _{\Xx^{\rm an}}(F, |A^{\rm an}_{\cdot}|).
$$
If $G= |A_{\cdot}|$ then by the above discussion $G^{\rm an}$ is
Illusie weak-equivalent to $|A^{\rm an}_{\cdot}|$ so by Theorem
\ref{main} we get
$$
\Re _{\Xx}(F, G) \cong
\Re _{\Xx^{\rm an}}(F, G^{\rm an})
$$
as claimed.
\eop

\numero{Examples}

\subnumero{Quotient stacks}

Consider the case of the site $\Xx$ of schemes over $\cc$.
Suppose $X$ is a scheme and $H$ is an algebraic group acting on
$X$.  Let $G:= X/H$ be the quotient pre-stack, in other words
$G(Y)$ is the groupoid corresponding to the action of $H(Y)$ on
$X(Y)$. Let $G'$ be the quotient stack, which is the stack
associated to the pre-stack $G$, i.e. the homotopy-sheafification.
The morphism $G\rightarrow G'$ is an Illusie weak equivalence.
We are interested in calculating
$$
(X/H)^{\rm top}:= \Re _{\Xx}(F,G')
$$
where $F$ is the functor $Y\mapsto Y^{\rm top}$.
By Theorem \ref{main} it suffices to calculate
$\Re _{\Xx}(F,G')$.  Note that $G'$ is the object-by-object
realization of the simplicial presheaf
$$
A_n := X\times H \times \ldots \times H.
$$
We  get that $(X/H)^{\rm top}$ is the realization
of the simplicial space
$$
A_n ^{\rm top}= X^{\rm top}\times H^{\rm top}
\times \ldots \times H^{\rm top}.
$$
Let $EH^{\rm top}$ denote a contractible space on which $H^{\rm
top}$ acts topologically freely, so that the quotient is $BH^{\rm
top}$. Then put
$$
S_n:=A_n\times EH^{\rm top}
$$
with simplicial structure modified to reflect the diagonal action
of $H^{\rm top}$ on $X^{\rm top}\times EH^{\rm top}$ (this action
is topologically free). We have a morphism of simplicial spaces
$S_n\rightarrow A_n^{\rm top}$ which is a weak equivalence at each
stage; thus it gives a weak equivalence of realizations. But the
realization of  $S_{\cdot}$ is just the quotient $(X^{\rm top}\times
EH^{\rm top})/H^{\rm top}$. The conclusion is that there is
a weak homotopy equivalence
$$
(X/H)^{\rm top}\cong (X^{\rm top}\times
EH^{\rm top})/H^{\rm top}.
$$
Not very surprising, after all...

\subnumero{Gerbs}

Suppose $Z$ is a connected scheme and $Y\rightarrow Z$ is a morphism
of algebraic stacks.
Recall that $Y$
is a {\em gerb} or {\em gerb over $Z$} if there is an
etale covering $U\rightarrow Z$ and a smooth group scheme $H_U$
over $U$ such that $Y\times _ZU$ is equivalent to the stack
associated to the pre-stack $K(H_U/U, 1)$ (this notation
means the stack which to a point $X\rightarrow U$ associates the
groupoid with one object and automorphism group $H_U(X)$).  Note
that in this case $H_U$ descends to a section $\Hh $ of the sheaf of
sets of isomorphism classes of group schemes over $Z$ (which we can
think of as being a group scheme over $Z$ but which is only defined
up to isomorphism---and this etale locally).

We say that $Y$ is a {\em split gerb} if there is a group scheme
$H$ over $Z$ and if $Y$ is the stack associated to the pre-stack
$K(H/Z,1)$ over $Z$.

\begin{lemma}
\label{splitGerb}
If $Z$ is a connected scheme of finite type  and
$Y\rightarrow Z$ is a split gerb with
group scheme $H$ which is smooth and a fibration over $Z$, then $Y^{\rm
top}\rightarrow Z^{\rm top}$ is homotopic to a fibration with fiber $B(H_z^{\rm
top})$ where $H_z$ is the fiber of $H$ over a basepoint $z\in Z$ and
$B(H_z^{\rm top})$ denotes the classifying space of the topological
group $H_z$.
\end{lemma}
{\em Proof:}
Express $Y$ as the realization of the simplicial set
$$
A_n = H\times _Z \ldots \times _ZH.
$$
Then $Y^{\rm top}$ is equivalent to the realization $|A^{\rm top}_n|$
where here, since the components of $A_n$ are schemes, we can just take the
usual topological realization at each stage. This realization is a fibration
(since  $H$ is a a fibration over $Z$) and the fiber is just the standard
expression for $B(H_z^{\rm top})$.
\eop

{\em Caution:}  If we take $Z$ equal to the affine line and let $H'$ be the
constant group scheme with fiber $GL(n)$ for example, then let $H$ be
obtained by blowing up the identity element over the origin in the affine line
and throwing out the part at infinity, we obtain a smooth group scheme over the
affine line whose generic point is $GL(n)$ and whose special point is $Lie
(GL(n))$ considered as an abelian group. In this case the topological type of
the fiber changes so our realization will no  longer be a fibration over
$Z^{\rm top}$. This is the reason for the condition that $H\rightarrow Z$
should be a fibration, in the hypothesis of the lemma.

Go back to the case of a general gerb $Y\rightarrow Z$ and let
$U\rightarrow Z$ be the etale covering over which it splits, given by the
definition. Assume as before that the group scheme type $\Hh$ over $Z$ is
smooth and a fibration (i.e. the group scheme $H_U$ is smooth and a fibration
over $U$).

We obtain a simplicial stack  $$ A_n := Y\times _ZU\times _Z \ldots
\times _ZU
$$
(or equivalently a simplicial presheaf of spaces).  The original
stack $Y$ is Illusie-equivalent to the object-by-object realization
$|A_{\cdot}|$. By \ref{simplicial} we can calculate the global realization of
$|A_{\cdot}|$  by taking the global realization at each stage
and then realizing the resulting simplicial space,
$$
|A_{\cdot}|^{\rm top} =|A_{\cdot}^{\rm top}|.
$$
By the main
theorem \ref{main} we get
$$
Y^{\rm top} \cong |A_{\cdot}^{\rm top}|.
$$
On the other hand let $C_n:= U\times _Z\ldots \times _ZU$.
It is a simplicial scheme realizing $Z$, so again applying Theorem
\ref{main} and the previous discussion we have
$$
|C_{\cdot}^{\rm top}|\stackrel{\cong}{\rightarrow} Z^{\rm top}.
$$
Via these equivalences, the morphism $Y^{\rm top}\rightarrow Z^{\rm
top}$ is homotopic to the morphism
$$
|A_{\cdot}^{\rm top}|\rightarrow |C_{\cdot}^{\rm top}|
$$
induced by the morphism of simplicial stacks $A_{\cdot}\rightarrow
C_{\cdot}$.

Note that there is a sheaf of groups $H_{C_n}$ on $C_n$, the
pullback of $H_U$; and we have that $A_n$ is the stack associated
to the pre-stack $K(H_{C_n}/C_n, 1)$.  In particular there is a
morphism (Illusie equivalence) $K(H_{C_n}/C_n, 1)\rightarrow A_n$
over $C_n$. By Theorem \ref{main} this induces a weak homotopy
equivalence
$$
K(H_{C_n}/C_n, 1)^{\rm top}\rightarrow A_n^{\rm top}
$$
compatible with the morphism to $C_n^{\rm top}$.
From our treatment of the case of a split gerb, the map
$A_n^{\rm top}\rightarrow C_n^{\rm top}$ is a fibration with fiber
$B(H_z^{\rm top})$. Note that we can choose the basepoint $z\in Z$ to
lift to any connected component of $U$, and  the fibers of $H_{C_n}$
over lifts  of $z$ are all isomorphic---so we get a
fibration with the same fiber over all connected components.
Note also that the morphisms in the simplicial space
$C_{\cdot}^{\rm top}$ induce isomorphisms of the fibers.   This in
turn implies that the realization
$|A_{\cdot}^{\rm top}|$ is a fibration over  $|C_{\cdot}^{\rm top}|$
with fiber $B(H_z^{\rm top})$.  We have proved the following
statement.

\begin{proposition}
\label{gerb}
Suppose $Z$ is a connected scheme.
If $Y\rightarrow Z$ is a gerb with
isomorphism-type of group schemes
$\Hh$ smooth and a fibration over  $Z$, let $H_z$ be a representative for the
fiber of a group scheme of type $\Hh$ over a basepoint $z\in Z$.  Then
$Y^{\rm top}\rightarrow Z^{\rm top}$ is a fibration with fiber
$B(H_z^{\rm top})$.
\end{proposition}
\eop

Look now at the case of Deligne-Mumford gerbs (i.e. gerbs which
are Deligne-Mumford stacks).  Note that a gerb $Y\rightarrow Z$
is
Deligne-Mumford if and only if the group-scheme type $\Hh$ is that
of a finite group $H$.   In this case $H^{\rm top}=H$ and we get
that $Y^{\rm top} \rightarrow Z^{\rm top}$ is a fibration with fiber
$B(H)$.    We have the following converse saying that to give
the gerb is the same as to give the fibration. Some notation: $Aut
(BH)$ denotes the $H$-space (or actually $A_{\infty}$-space) of
self-homotopy equivalences of $BH$, and $B(Aut(BH))$ denotes its
delooping. Recall that
$$
\pi _0(Aut (BH))= Out (H)
$$
is the group of outer automorphisms of $H$ and
$$
\pi _1(Aut (BH))=Z(H)
$$
is the center of $H$.  Thus $B(Aut(BH))$ is a $2$-truncated
connected space with $\pi _1=Out (H)$ and $\pi _2=Z(H)$ (it is
associated to the standard example of a crossed-module $H\rightarrow
Aut(H)$).

\begin{proposition}
\label{converse}
Suppose $Z$ is a connected scheme and $H$ is a finite group.
Then the $2$-category of Deligne-Mumford gerbs $Y\rightarrow Z$ of type $H$ is
equivalent to the $2$-category of fibrations over $Z^{\rm top}$ with fiber
$H$ or
equivalently the Poincar\'e $2$-category of the space $Hom (Z^{\rm top}, B(Aut
(BH))$.   \end{proposition}
{\em Proof:}
This is a comparison theorem between etale cohomology and usual cohomology with
coefficients in the $2$-category $B(Aut (BH))$. Using  Postnikov devissage on
this last $2$-category it suffices to prove the comparison theorem for the
cohomology with coefficients in $\pi _1(B(Aut(BH)))$ and $\pi _2(B(Aut (BH)))$.
But these are finite groups
so the
comparison theorem between etale and usual cohomology holds.
\eop

\subnumero{Deligne-Mumford stacks}
We show how to use  the
comparison theorem \ref{comparison} and Theorem \ref{main} applied
to $\Xx ^{\rm an}$ to give a cut-and-paste description of the
topological realization of a Deligne-Mumford stack.

Suppose $Y$ is a Deligne-Mumford stack.  We have the following result:
\begin{proposition}
There exists a {\em coarse moduli space} $Z$ for $Y$, in other words
a scheme with a morphism $Y\rightarrow Z$ which is universal for
morphisms from $Y$ to schemes. Furthermore, locally in the etale
topology of $Z$ we can express $Y$ as a quotient stack by a finite
group action and $Z$ as the corresponding quotient scheme.
\end{proposition}
Mochizuki alludes to this result in \cite{Mochizuki} and gives
\cite{Faltings-Chai} as reference.
\eop

{\em Assumption:}  We assume that $Z$ can be stratified as a disjoint
union of connected smooth locally closed subschemes
$$
Z= \bigcup Z_{\beta}
$$
such that the stabilizer subgroups of the group action are
locally constant over $Z_{\beta}$ and $Z$ (as stratified) is
equisingular along the $Z_{\beta}$.

This assumption is probably always true, but we don't try to prove it
here since it is not our purpose to get into a long discussion of
stratifications.

Let $Y_{\beta}$ denote the
inverse image of $Z_{\beta}$ in $Y$, and let $H_{\beta}$ denote the
isomorphism class of stabilizer group over $Z_{\beta}$.
Note that $Y_{\beta}\rightarrow Z_{\beta}$ is
a gerb with group $H_{\beta}$.

Let $Z_{\dim \leq n}$ denote the union of strata of dimension $\leq n$.
It is a closed subset, and
$$
Z_{\dim = n}:= Z_{\dim \leq n} - Z_{\dim \leq n-1}
$$
is a disjoint union of the strata of dimension $n$. Use the same
notations in $Y$.

The topological realization $Y^{\rm top}_{\dim = n}$ is a fibration
over $Z^{\rm top}_{\dim = n}$ whose fiber over a connected component
$Z_{\beta}^{\rm top}$ is $K(H_{\beta}, 1)$.  The gerb
$Y_{\dim = n}$ is determined by $Z_{\dim = n}$ and this fibration.

Let $T_{\leq n}$ be a good (as usual in the theory of singular spaces
and stratifications) tubular neighborhood of
$Z^{\rm an}_{\dim \leq n}$ in $Z^{\rm an}$.  Let $V_{\leq n}$ be the
inverse image of $T_n$ in $Y^{\rm an}$.

Let $T_{=n}$ be the complement of $T_{\leq n-1}$ in $T_{\leq n}$ (and
we assume that this is nicely arranged so as to be a tubular
neigbhorhood of a subset which is essentially $Z_{\dim =n}$). Again
let $V_{=n}$ be the inverse image in the stack.

{\em Inductive claim:} The retraction of $T_{\leq n}^{\rm top}$ to
$Z_{\dim
\leq n}^{\rm top}$ lifts homotopically to a retraction from
$V_{\leq n}^{\rm top}$ to $Y_{\dim
\leq n}^{\rm top}$. This retraction preserves $V_{=n}^{\rm top}$ over
$T_{=n}^{\rm top}$.

In particular $V_{=n}^{\rm top}$ is homotopy equivalent to
$Y_{\dim =
n}^{\rm top}$, which as we have said above is a fibration over
$Z_{\dim =n}^{\rm top}$ with fibers $K(H_{\beta}, 1)$.

Applying the main theorem to the presentation of the stack $V_{
\leq n}$ corresponding to the covering by $V_{\leq n-1}$ and
$V_{=n}$  we find that $V_{\leq n}$ is obtained by
glueing $V_{=n}^{\rm top}$ (which we understand) to $V_{\leq n}^{\rm
top}$ (which we suppose we understand by induction) along a boundary
which is homotopic to the inverse image in
$Y$ of the boundary of $T_{\leq n-1}\cap Z_{\dim
\leq n}$. This boundary piece is again a gerb.

By looking at this closely one can arrange to have the
inductive claim for $n$ once it is known for $n-1$.

Finally when we get to $n= \dim (Z)$ we are done and we have
constructed $Y^{\rm top}$ by a sequence of ``cut and paste''
operations. To recapitulate this sequence of operations, we get
$Y_{\leq n}^{\rm top}$ from $Y_{\leq n-1}^{\rm top}$ by adding on a
space which is a fibration over $Z_{\dim = n}$ (the union of open
strata of dimension $n$) with fibers of the form $K(H_{\beta}, 1)$ for
various finite groups $H_{\beta}$.

{\em Remark:}  The whole of the above discussion goes through
equally well for any algebraic stack $Y$, if we assume the
existence of a morphism $Y\rightarrow Z$, and assuming the existence
of a nice stratification as in the assumption stated above (in
particular this implies that $Y\rightarrow Z$ gives an isomorphism
of $Spec (\cc )$-valued points).  Note that this assumption is not
automatic, since for example it doesn't hold for
the moduli stacks of vector bundles and the like (for in those
cases the coarse moduli space is no longer an
isomorphism on points).

{\em Problem:} Give a ``cut and paste''
description of the topological realization of the moduli stacks of
vector bundles (as opposed to the quotient space description
arising from the expression of the moduli stack as a quotient of
the Hilbert scheme).

{\em 3. The case of curves}
Finally we treat very explicitly what happens in the case of curves.
Suppose $Y$ is a smooth Deligne-Mumford stack of dimension $1$.
Let $\pi :Y\rightarrow Z$ be the coarse moduli space. Since $Y$ is
normal, the coarse moduli space is also normal so $Z$ is a smooth
curve.  Let $P_i\in Z$ be the points of ramification of $\pi $ and
let $n_i$ be the ramification indices.  Define the intermediate
stack $Y\rightarrow W \rightarrow Z$ to be the orbicurve with
ramification $n_i$ at $P_i$ over $Z$.  It is given by local charts
which ramify at $P_i$ with degree $n_i$ and are etale elsewhere
(but don't meet the $P_j$) $j\neq i$).  The morphism $Y\rightarrow
W$ is a gerb with group $H$ where $\pi ^{-1}(x)=K(H,1)$ for a
general point $x\in Z$.

Applying the cut-and-paste construction above, one can see that
$W^{\rm top}$ is obtained by glueing the lens spaces $K(\zz /n_i \zz
, 1)$ into the space $(W-\{ P_i\} )^{\rm top}$ at the punctures
$P_i$, with the loop around the puncture going to the standard
generator for the fundamental group of  the lens space.
Then from the discussion of gerbs, $Y^{\rm top} \rightarrow W^{\rm
top}$ is a fibration with fiber $K(H, 1)$.  Lemma \ref{converse}
says that the classifying element $\eta$ for this fibration in the
appropriate classifying space $Hom (W^{\rm top}, B\, Aut (K(H,1)))$
gives exactly the data of the gerb $Y\rightarrow W$.  In particular,
the data of the stack $Y$ can be given by the curve $Z$, the points
$P_i$ and ramification indices $n_i$ (which results up to now in a
space $W^{\rm top}$), the group $H$ (up to isomorphism) and the
map $\eta :  W^{\rm top}\rightarrow B\, Aut (K(H, 1))$.
The map $\eta$ yields first a map $\eta _1:\pi _1(W^{\rm
top},x)\rightarrow  Out (H)$ describing the twisting of the group
$H$ as one moves around in $W$, and second, an element $\eta
_2$ of $H^2(W^{\rm top}, Z(H))$ where the coefficient system (the
center of $H$) is twisted by $\eta_1$.

\numero{The generalisation to $n$-topoi}

This section is a direct result of conversations (email) with K.
Behrend (in the course of preparations for the Trento stack
school) and C. Telemann, who is thinking about this type of thing
in connection with $G$-bundles on a curve. At the current writing
the arguments below are only sketches---I haven't checked the details.

For this section we make a major change of notation: now $\Delta$ will denote
the simplicial category whose objects are ordered sets of the form $[n]=\{
0,\ldots , n\}$ and whose morphisms are  increasing morphisms of ordered sets.
Thus $\Delta ^n$ now denotes the cartesian product of the category $\Delta$
with itself $n$ times.

Recall that an $n$-category is a functor $\alpha : \Delta
^n\rightarrow Sets$ which satisfies certain conditions
\cite{Tamsamani}. The set of objects of $\alpha$ is $\alpha _{0,
\ldots , 0}$.  Let $\alpha _m (x_0,\ldots , x_m)$ denote the
$n-1$-category of composable $m$-tuples with objects $X_0,\ldots ,
x_m$.

We assume known a definition of the $n$-category of functors
$Hom (\alpha , \beta )$ between two $n$-categories, such that the
composition  $$ Hom (\alpha , \beta )\times Hom (\beta , \gamma
)\rightarrow Hom (\alpha , \gamma ) $$
is strictly associative.  This seems to exist from some
preliminary thoughts on the subject, although noone has yet checked
the details.

If $\alpha$ is an $n$-category we can define a new $n$-category
$\alpha \ast$ as the $n$-category with one final object attached.
The objects of $\alpha \ast$ are those of $\alpha$ plus one more
object denoted $e$, and
$$
\alpha \ast _n(x_0,\ldots , x_{j} , e,\ldots , e):=
\alpha _{j}(x_0,\ldots , x_{j}),
$$
and $\alpha \ast _n(x_0,\ldots , x_n)=\emptyset $ if $x_i=e$ and
$x_j\neq e$ for $j>i$---this gives the definition of the
$n$-category $\alpha \ast$.

Similarly we can define the $n$-category $\ast \alpha$ by attaching
one  initial object.

These satisfy homotopy-universal properties which we leave to be
elucidated later.  Heuristically these properties will say that if
$f\alpha \rightarrow \beta$ is a functor of $n$-categories and if
$u\in Ob (\beta )$ is an object with a natural transformation
$f\rightarrow cu$ (where $cu: \alpha \rightarrow \beta$ is the
constant functor with values $u$) then there is an essentially
unique functor $\alpha \ast \rightarrow \beta $ sending $e$ to $u$
(or to an object equivalent to $u$ provided with the equivalence).
The same goes with arrows reversed for $\ast \alpha$.

In view of these universal properties it is reasonable to say that
{\em $\beta$ admits direct limits indexed by $\alpha$} if the
functor of $n$-categories (restriction)
$$
Hom (\alpha \ast , \beta )\rightarrow Hom (\alpha , \beta )
$$
is an equivalence.  Similarly we say that {\em $\beta$ admits
inverse limits indexed by $\alpha$} if the functor
$$
Hom (\ast \alpha , \beta )\rightarrow Hom (\alpha , \beta )
$$
is an equivalence.

Needless to say, the limits we are talking about
here correspond to what is sometimes called ``holim'' in the usual
topological situation.
The category $Top$ admits a natural structure of $n$-category for
any $n$ (denoted $Top_n$); the objects are the topological spaces
and the morphism $n-1$-category from $X$ to $Y$ is defined to be the
Poincar\'e $n-1$ category of the space $Hom^{\rm top}(X,Y)$. Note that
this structure of $n$-category is what might usually be called the
$n$-category of $n-1$-truncated spaces.  It is equivalent to the
$n$-category of $n-1$-groupoids. The usual $holim$ (direct and
inverse) are limits in the above sense, indexed over
$1$-categories, in the $n$-category $Top_n$.

We adopt the working definition\footnote{Of course this
circumlocution means that we haven't yet checked the details!}
that an {\em $n$-topos} is an $n$-category which admits direct and
inverse limits indexed by $n$-categories.

The  $n$-category $Top _n$ is the basic example of an $n$-topos,
extending the topos of sets (which is the case $n=1$). More
generally the $n$-category of $n-1$-stacks of groupoids over a
site $\Xx$ is an $n$-topos which we call the {\em $n$-topos
associated to $\Xx$}.

If $\tau$ is an $n$-topos the existence of limits automatically
gives a functor $Top _n\rightarrow \tau$. In the case of the
$n$-topos of $n-1$-stacks over $\Xx$ this is the functor taking an
$n-1$-groupoid to the associated constant stack.

\subnumero{Realizations with values in an $n$-topos}
We now recast our discussion of topological realization in
terms of functors with values in an $n$-topos $\tau$.
Suppose $\Xx$ is a category and $F: \Xx \rightarrow \tau$ is a
covariant functor.  Suppose that $G$ is an $m$-stack.  Then we
obtain $\Re _{\Xx} (F,G)\in \tau$, functorial (in the
homotopic sense) in the two variables $F$ and $G$.  To define this,
we take essentially the same definition as before, but noting that
the explicit space that is constructed there is really just a
direct holim.  More precisely let $Fl (\Xx )$ denote the category
of morphisms in $\Xx$, we have
$$
Fl (\Xx )\rightarrow \Xx \times \Xx ^o
$$
but on the other hand $(F,G)$ (and the canonical $Top _m
\rightarrow \tau$) gives a functor $\Xx \times \Xx ^o \rightarrow
\tau$. Composing we obtain the functor $Fl (\Xx )\rightarrow \tau$
and the realization is its direct limit in $\tau$.  The morphisms
which we could construct explicitly before are now only canonically
defined up to canonically defined homotopy up to etc (again the
details still have to be worked out carefully).

Suppose $\Xx$ is a site. As before we say that $F$ satisfies {\em
covariant descent} if for any $X\in \Xx$ and sieve $\Bb \subset
\Xx /X$ the morphism $\Re _{\Bb} (F|_{\Bb}, \ast _{\Bb})\rightarrow
F(X)$ is an equivalence in $\tau$.  The same proof as above goes
through to give the following.

\begin{theorem}
\label{topoi}
If $\tau$ is an $n$-topos, $\Xx$ a site, and $F: \Xx \rightarrow
\tau$ a covariant functor satisfying covariant descent and $G,G': \Xx
\rightarrow Top_m$ two presheaves of $m-1$-truncated spaces
with a morphism $G\rightarrow G'$ then the induced morphism
$\Re _{\Xx} (F, G)\rightarrow \Re _{\Xx} (F, G')$ is an equivalence
in $\tau$.  Consequently the functor $F$ extends to a functor from
the topos of $m-1$-stacks on $\Xx$ to $\tau$ commuting with
(homotopy!) direct limits.
\end{theorem}
\eop

In the case $\tau = Top _n$ we recover the previous Theorem
\ref{main}.

The following corollary extends to the $n$-topos of $n-1$-stacks the
usual definition of pullback of sheaves for a morphism of sites.

\begin{corollary}
\label{pullback}
Suppose $\Xx$ and $\Yy$ are sites and $F: \Xx \rightarrow \Yy$ is
a functor (which is sometimes called a morphism of sites from $\Yy$
to $\Xx$).  Then $F$ extends uniquely to a functor from
the $n$-topos
associated to $\Xx$ to the $n$-topos associated to $\Yy$,
compatible with homotopy direct limits.
\end{corollary}
\eop

{\em Question:} Where do inverse limits come into things?

{\em Remark:} In the case where $\tau$ is the $n$-topos of
$n-1$-truncated presheaves of spaces over a site $\Yy$ (such as for
Corollary \ref{pullback}) one can again do everything by hand as in
the first part of the paper (the details are pretty much identical)
so in this case one doesn't need to worry about getting all of the
foundational material about $n$-categories straight.

C. Telemann has a nice interpretation of the realization obtained from the
functor $F: X\mapsto X^{\rm top}$.  He points out that on the analytic site
$\Xx ^{\rm an}$ the topological realization can be defined ``internally'' as in
the etale homotopy theory of Artin-Mazur, Friedlander.  In our notations I
think that this means that the constant functor $F_0(X):= \ast$ has a canonical
replacement by a functor  $F$ which satisfies covariant descent (one just
enforces covariant descent in a way analogous to the operation $\gamma$ for
enforcing the homotopy-sheaf condition).  This functor is equivalent to the
standard topological realization functor.  Then Telemann points out that by
pulling back a stack from the algebraic to the analytic site (with pullback
defined as in Corollary \ref{pullback}) and applying this internal topological
realization functor we obtain our topological realization.

{\em Remark:}  The operation $G\mapsto G^{\rm an}$ that we have defined
slightly differently should be the same as the pullback via the morphism of
sites
$$
\mbox{(analytic)} \rightarrow \mbox{(algebraic)}
$$
(which corresponds to the standard functor $\Xx \rightarrow \Xx ^{\rm an}$).

\end{document}